\documentclass[12pt,twoside]{article} 
\usepackage{amssymb,mathrsfs,graphicx}

\newtheorem{theorem}{Theorem}[section]
\newtheorem{lemma}[theorem]{Lemma}

\newtheorem{definition}[theorem]{Definition}

\newtheorem{rmrk}[theorem]{Remark}

\makeatletter
\@addtoreset{equation}{section}
\makeatother

\newenvironment{remark}
{\begin{rmrk} \em}
{\end{rmrk}}

\newcommand{\fn} {function}
\newcommand{\me} {measure}
\newcommand{\tr} {trajector}
\newcommand{\erg} {ergodic}
\newcommand{\sy} {system}
\newcommand{\hyp} {hyperbolic}
\newcommand{\pr} {probability}
\newcommand{\ra} {random}
\newcommand{\dsy} {dynamical system}

\newcommand{\R} {\mathbb{R}}

\newcommand{\Z} {\mathbb{Z}}
\newcommand{\N} {\mathbb{N}}
\newcommand{\qed} {\hfill {\small Q.E.D.} \par\medskip}
\newcommand{\skippar} {\par\medskip}

\newcommand{\proof} {\noindent \textsc{Proof.} }

\newcommand{\article}[3] {\textsc{{#1}}, {\itshape {#2}}, {{#3}}.}
\newcommand{\book}[3] {\textsc{{#1}}, {\itshape {#2}}, {{#3}}.}
\newcommand{\vol} {\textbf}
\newcommand{\eps} {\varepsilon}
\newcommand{\rset}[2] {\left\{ #1 \: \left| \: #2 \right. \! \right\} }
\newcommand{\lset}[2] {\left\{ \left. \! #1 \: \right| \: #2 \right\} }

\renewcommand{\iff} {if and only if\ }

\newcommand{\rw} {random walk}
\newcommand{\en} {environment}

\newcommand{\om} {\omega}

\newcommand{\ps} {\mathcal{M}}
\newcommand{\ma} {T}

\begin{document}

\title{\textbf{Central Limit Theorem and recurrence \\ 
for random walks in \\ 
bistochastic random environments}}

\author{\textsc{Marco Lenci}
\thanks{
Dipartimento di Matematica, 
Universit\`a di Bologna, 
Piazza di Porta S.\ Donato 5, 
40126 Bologna, Italy.
E-mail: \texttt{lenci@dm.unibo.it} } 
}

\date{October 2008}

\maketitle

\begin{abstract}
  We prove the annealed Central Limit Theorem for random walks in
  bi\-stochastic random environments on $\Z^d$ with zero local
  drift. The proof is based on a ``dynamicist's interpretation'' of the
  system, and requires a much weaker condition than the customary
  uniform ellipticity. Moreover, recurrence is derived for $d \le 2$.

\end{abstract}

\section{Introduction}
\label{sec-intro}

We study the Central Limit Theorem (CLT) and the recurrence properties
of a certain class of \rw s in \ra\ \en s (RWREs), namely the \rw s in
bistochastic \en s with zero local drift.  Although this class is
fairly general, it is not as general as we can prove theorems for, and
certainly not new, having been previously investigated at least by
Koslov \cite{k} and Komorowski and Olla \cite{ko}.

In fact, the purpose of this note is not to give our most original
results for the amplest class of RWREs (we will take this point of
view in another paper \cite{l3}; see also Section \ref{sec-gen}
below), but rather to present a technique, at the frontier of \dsy s
and \pr\ theory, that can deliver known results more easily than the
current methods, and often improve on them. Perhaps more importantly,
it provides a unifying view on the relation between diffusive behavior
and recurrence for \ra\ as well as deterministic dynamics \cite{l1,
l2}.  Our method hinges in part on a beautiful theorem by Schmidt
\cite{s} (a somewhat weaker version of which has been independently
achieved by Conze \cite{c}) on the recurrence of commutative cocycles
over \erg\ \dsy s.

Let us describe the \sy: We deal with \rw s on $\Z^d$ defined by a
bistochastic matrix $p := \{ p_{xy} \}_{x,y \in \Z^d}$ of transition
probabilites. This means that
\begin{equation}
  \label{bi-en}
  \forall y \in \Z^d, \qquad \sum_{x \in \Z^d} p_{xy} =1,
\end{equation}
together of course with the trasposed condition, customarily called
\emph{normalization}, which ensures that, $\forall x \in \Z^d$, $y
\mapsto p_{xy}$ is indeed a \pr\ distribution on $\Z^d$.  For the sake
of simplicity we assume that there exists a finite $\Lambda \subset
\Z^d$ such that
\begin{equation}
  \label{lambda}
  p_{xy} = 0 \qquad \mbox{ if } y-x \not\in \Lambda.
\end{equation}
Our \rw s will have \emph{zero local drift}, in the sense that
\begin{equation}
  \label{zero-l-drift}
  \forall x \in \Z^d, \qquad \sum_{y \in \Z^d} (y-x) \, p_{xy} = 0.
\end{equation}

We assume that the transition matrix $p$ is itself \ra: $p = p(\om)$,
where $\om$ ranges in the \pr\ space $(\Omega, \Pi)$. We will
liberally call both $p(\om)$ and $\om$ the \emph{\ra\ \en}, or simply
the \en. It is not important what $\Omega$ actually is (although the
interested reader may look at Section \ref{subs-ex} of the Appendix)
but we make the fundamental hypothesis that it is acted upon by the
group $\{ \tau_z \}_{z \in \Z^d}$ of automorphisms w.r.t.\ $\Pi$. This
action is such that
\begin{equation}
  \label{zd-action}
  p_{xy}(\tau_z \om) = p_{x+z,y+z}(\om)
\end{equation}
and it is \erg\ (which is the minimal assumption for the \ra\ law on
the \en\ to have something to do with the dynamics of the
walker). Because of this, it is no loss of generality to require that
the walk always starts at $0$.

Our last hypothesis is the \emph{almost sure irreducibility} of $p$:
For $\Pi$-a.e.\ $\om \in \Omega$, for every $y \in \Z^d$, there exists
$n = n(\om,y)$ such that
\begin{equation}
  \label{as-irr}
  p^{(n)}_{0y} := \sum_{x_1, \ldots x_{n-1}} p_{0x_1} \, p_{x_1 x_2}
  \cdots\, p_{x_{n-1} y} > 0.
\end{equation}
(It is easy to see, via (\ref{zd-action}), that (\ref{as-irr})
guarantees mutual accessibility of any two points $x,y$ of $\Z^d$,
at least for a.a.\ \en s, whence the name `irreducibility'.) An
example of $(\Omega, \Pi)$ is given in Section \ref{subs-ex} of the
Appendix. 

\begin{remark}
  Notice how condition (\ref{as-irr}) is much weaker than the
  \emph{uniform ellipticity} assumed in most results on RWREs, namely
  the existence of a constant $\eps>0$ such that, for a.e.~$\om$,
  \begin{equation}
    \label{unif-ell}
    p_{0e} \ge \eps, \qquad \forall e \in \Z^d, |e|=1.
  \end{equation}
  To the author's knowledge, within the scope of the diffusive or
  recurrence properties of RWREs \cite{aks, la, ko}, only the papers
  by Berger and Biskup \cite{be, bb} do not (and cannot) require
  uniform ellipticity (cf.\ also \cite{z1, z2}).  Moreover,
  (\ref{as-irr}) can be further relaxed if stronger \erg\ properties
  hold for $(\Omega, \Pi, \{ \tau_z \})$, cf.\ Section \ref{sec-gen}.
\end{remark}

For the sake of mathematical rigor, we now give a formal definition of
our \sy\ in terms of standard objects of \pr\ theory, warning the
reader that the following construction is not the one that we will
work with in the rest of the paper. Indeed, from Section
\ref{sec-povp} onward, we will represent the above RWRE in terms of a
suitable \me-preserving \dsy.

At any rate, fixed $\om \in \Omega$, the \rw\ in the \en\ $\om$ is the
Markov chain $P_\om$ on $\Z^d$ defined by
\begin{eqnarray}
  && P_\om (X_0 = 0) = 1; \\
  && P_\om (X_{n+1} = y \,|\, X_n = x) = p_{xy}(\om). 
\end{eqnarray}
We take into account the complete \ra ness of the problem by studying
the stochastic process $\{ X_n \}_{n\in\N}$ w.r.t.\ the \emph{annealed
law}, , which is defined on $\Omega \times (\Z^d)^\N$ via
\begin{equation}
  \label{bbp}
  \mathbb{P} (B \times E) := \int_B \Pi(d\om) \, P_\om(E),
\end{equation}
where $B$ is a Borel set of $\Omega$ and $E$ a Borel set of
$(\Z^d)^\N$ (the latter being the space of the \tr ies, where $P_\om$
is defined).

The paper's main results are the annealed CLT, i.e., the CLT relative
to $\mathbb{P}$ (Theorem \ref{thm-clt}), and the almost sure
recurrence in dimension $d \le 2$, namely, the property that the \rw\
is recurrent in a.a.\ \en s (Theorem \ref{thm-rec}).

\skippar

The exposition is organized as follows: In Section \ref{sec-povp} we
introduce the \dsy\ that we use to represent our RWRE, whose \erg\
properties we study in Section \ref{sec-erg-clt}. In Section
\ref{sec-rec} we present Schmidt's result on recurrent cocycles and
apply it to our \sy. Finally, in Section \ref{sec-gen}, we draw some
brief conclusions about the present work and discuss how to generalize
it in a number of ways.

\bigskip
\noindent
\textbf{Acknowledgments.}\ I wish to thank Pierluigi Contucci and
Cristian Giardin\`a for inviting me to speak at the YEP-V Meeting and
for causing me to organize my ideas by requesting my contribution to
this Special Issue of JMP.

\section{The point of view of the particle}
\label{sec-povp}

Let us enumerate the elements of $\Lambda$, cf.\ (\ref{lambda}), as
$d_1, d_2, \ldots, d_N$.  Now let us fix $\om \in \Omega$. For $i = 1,
\ldots, N$, we define
\begin{eqnarray}
  \label{qi}
  && q_i = q_i(\om) := p_{0d_i}(\om) \\ 
  \label{qip}
  && q'_i = q'_i(\om) := p_{-d_i 0}(\om) = q_i (\tau_{-d_i} \om),
\end{eqnarray}
the last equality coming from (\ref{zd-action}).  By
(\ref{bi-en})-(\ref{lambda}),
\begin{equation}
  \label{sum-q}
  \sum_{i=1}^N q_i = \sum_{i=1}^N q'_i = 1.
\end{equation}
We then set $a_0 := 0$ and, recursively for $i=1, \ldots, N$, 
\begin{eqnarray}
  \label{ai}
  a_i = a_i(\om) &:=& a_{i-1} + q_i \\
  \label{ii}
  I_i = I_i(\om) &:=& [a_{i-1} , a_i).
\end{eqnarray}
By the first of the (\ref{sum-q}), $\{ I_i \}$ is a partition of $I :=
[0,1)$.  For $(s,\om) \in I \times \Omega$, let $i(s,\om)$ be the
unique $i$ such that $s \in I_i(\om)$. Setting
\begin{eqnarray}
  \label{def-d}
  D(s,\om) &:=& d_{i(s,\om)}, \\
  \label{def-phi}
  \phi(s,\om) &:=& q_{i(s,\om)}^{-1} \left(s - a_{i(s,\om)} \right),
\end{eqnarray}
defines the \fn s $D: I \times \Omega \longrightarrow \Lambda$ and
$\phi: I \times \Omega \longrightarrow I$.  For reasons that will be
clear momentarily, $D$ is called the \emph{displacement \fn} and
$\phi$ is called the \emph{internal dynamics}, or the \emph{map on the
fibers}. It is apparent that $\phi(\cdot, \om)$ is a
piecewise-linear, at most $N$-to-1 map of $I$ onto itself. More
precisely, it is the perfect Markov map $I \longrightarrow I$ relative
to the partition $\{ I_i(\om) \}$.

The \dsy\ we study for the rest of the paper is the triple $(\ps, \mu,
\ma)$, where $\ps := I \times \Omega$, $\mu := m \times \Pi$ (having
denoted by $m$ the Lebesgue \me\ on $I$), and $\ma : \ps
\longrightarrow \ps$ is given by
\begin{equation}
  \label{def-ma}
  \ma (s,\om) := \left( \phi(s,\om) , \tau_{D(s,\om)} (\om) \right).
\end{equation}
This \sy\ is called the \emph{point-of-view-of-the-particle \dsy}
and the reason can be explained as follows.

Fix $\om \in \Omega$ and a \ra\ $s \in I$ w.r.t.\ $m$. The \pr\ that
$s \in I_i(\om)$ is $m(I_i) = q_i$, which, in terms of our \rw, is
exactly the \pr\ that a particle placed in the origin of $\Z^d$,
endowed with the \en\ $p(\om)$, jumps by a quantity $d_i$. Then, back
to the \dsy, condition the \me\ $m$ to $I_i$. Calling $(s_1,\om_1) :=
\ma (s,\om)$, we see that, upon conditioning, $s_1$ ranges in $I$ with
law $m$. Therefore, in a sense, the variable $s$ (which we may call
the \emph{internal variable}) has ``refreshed'' itself. Furthermore,
$\om_1$ is the translation of $\om$ in the opposite direction to $d_i
= D(s,\om)$, cf.\ (\ref{zd-action}). Hence we can imagine that we have
reset the \sy\ to a new initial condition $(s_1,\om_1)$, corresponding
to the particle sitting in $0 \in \Z^d$ and subject to the \en\
$p(\om_1)$. Applying the same reasoning to $(s_2,\om_2) := \ma
(s_1,\om_1)$, and so on, shows that we are following the motion of the
particle in the reference \sy\ of the particle itself, whence the
`point of view of the particle'.

In any case, it should be clear that the stochastic process $\{ X_n
\}$, with $X_0 := 0$ and, for $n \ge 1$,
\begin{equation}
  \label{def-xn-dk}
  X_n(s,\om) := \sum_{k=0}^{n-1} D_k(s,\om) := \sum_{k=0}^{n-1} D
  \circ \ma^k (s,\om),
\end{equation}
is precisely the \rw\ in the \en\ $p(\om)$. The definition of $\mu$
entails that $(\ps, \mu, \ma)$ describes all the realizations of the
\rw\ in \emph{all} the \en s, w.r.t.\ a \me\ that, in the language of
Section \ref{sec-intro}, is expressed by (\ref{bbp}). This is our
RWRE.

\section{Ergodic properties and Central Limit Theorem}
\label{sec-erg-clt}

In this section we study the stochastic properties of the \dsy\
defined above, starting with the most basic, the invariance of the
\me.

\begin{lemma}
  $T$ preserves $\mu$.
\end{lemma}

\proof Without loss of generality, it is sufficient to prove that
$\mu(\ma^{-1} A) = \mu(A)$ for sets of the type $A = [b,c] \times
B$, where $B$ is a measurable subset of $\Omega$. For added
simplicity, we may assume that every $(s,\om) \in A$ has exactly $N$
counterimages (the other cases can be considered degenerate versions
of the one we are assuming).

By direct inspection of the map (\ref{def-ma}), we can write $\ma^{-1}
A = \bigcup_{i=1}^N A'_i$, where
\begin{equation}
  \label{meas10}
  A'_i := \rset{ (s',\om') } {\om' \in \tau_{-d_i}(B), \, s' \in 
  [ a_i(\om') + q_i(\om') b \,,\, a_i(\om') + q_i(\om') c] },
\end{equation}
cf.\ (\ref{ai}), (\ref{def-phi}). These sets are pairwise disjoint
because, by construction, they belong to different level sets of the
\fn\ $D$. From (\ref{meas10}),
\begin{equation}
  \label{meas20}
  \mu(A'_i) = \int_{\tau_{-d_i}(B)} \! q_i(\om') (c-b) \, \Pi(d\om') = 
  (c-b) \int_B q'_i(\om) \Pi(d\om),
\end{equation}
having used, in the second equality, (\ref{qip}) and the
$\tau$-invariance of $\Pi$. Summing the above over $i=1, \ldots, N$,
with the help of the second of the (\ref{sum-q}), yields $(c-b) \Pi(B)
= \mu(A)$.  
\qed

Let us call \emph{horizontal fiber} of $\ps$ any segment of the type
$I_\om := I \times \{ \om \}$, and indicate by $m_\om$ the Lebesgue
\me\ on it. Also, given a positive integer $n$ and a multi-index
$\mathbf{i} := (i_0, i_1, ..., i_{n-1}) \in \mathcal{I}^n := \{ 1, 2,
\ldots, N \}^n$, we set
\begin{equation}
  I_{\om,\mathbf{i}} := \rset{(s,\om) \in \ps} {D_k (s,\om) =
  d_{i_k},\: \forall k=0, \ldots, n-1},
\end{equation}
where $D_k$ is defined in (\ref{def-xn-dk}). Finally, we denote by
$I_\mathbf{i} = I_\mathbf{i} (\omega)$ the interval of $I$
corresponding to $I_{\om,\mathbf{i}}$ via the natural isomorphism
$I_\om \longrightarrow I$.

It is easy to ascertain that $\{ I_\mathbf{i} \}_{\mathbf{i} \in
\mathcal{I}^n}$ partitions $I$ into $N^n$ intervals (some of which
may be empty), each corresponding to one of the realizations of the
\rw\ $\{ X_k \}_{k=0}^N$ in the \en\ $\om$, in such a way that $m
(I_\mathbf{i})$ is the \pr\ of the corresponding realization.  In
analogy with the previous notation, we indicate with
$I_{\mathbf{i}_n(s,\om)}$ the element of said partition that contains
$s$.

\begin{lemma}
  \label{lem-hyp}
  For a.a.\ $(s,\om) \in \ps$, $m( I_{\mathbf{i}_n(s,\om)} (\om))$
  vanishes exponentially fast, as $n \to \infty$.
\end{lemma}

\proof First of all, we introduce a notation that will be convenient
for this and other proofs. For $(s,\om) \in \ps$ and $k \in \N$, we
write
\begin{equation}
  \label{sk-omk}
  (s_k,\om_k) := \ma^k(s,\om).
\end{equation}
Now define 
\begin{equation}
  \label{erg5}
  f(s,\om) := \log q_{i(s,\om)}^{-1} (\om) = -\log m(I_{i(s,\om)}
  (\om)) 
\end{equation}
(with the convention that $\log 0 = -\infty$). Then $f(s,\om) \ge 0$,
the equality holding only when $\{ I_i(\om) \}$ is the trivial
partition of $I$, mod $m$, i.e., when $p_{0y}(\om) = \delta_{yy_0}$,
for some $y_0$. By (\ref{zero-l-drift}) it must be $y_0 = 0$. But this
can only happen for a negligible set of $\om$, due to the almost sure
irreducibility (\ref{as-irr}). 

Thus, $f > 0$ a.e. A well-known corollary of the Birkhoff Theorem
ensures that 
\begin{equation}
  \label{erg8}
  f^+ (s,\om) := \lim_{n \to \infty} \frac1n \sum_{k=0}^{n-1} f
  (s_k,\om_k) > 0
\end{equation}
as well, for a.a.\ $(s,\om) \in \ps$. On the other hand, from what we
have discussed earlier, it is easy to verify that, for $n \ge 1$,
\begin{equation}
  \label{erg10}
  m \left( I_{\mathbf{i}_n(s,\om)} (\om) \right) = \exp \left( -
  \sum_{k=0}^{n-1} f (s_k,\om_k) \right).
\end{equation}
The combination of (\ref{erg8}) and (\ref{erg10}) yields the
assertion.
\qed

\begin{lemma}
  \label{lem-erg}
  The \erg\ components of $(\ps, \mu, \ma)$ contain whole horizontal
  fibers, that is, every invariant set is of the form $I \times B$,
  mod $\mu$, where $B$ is a measurable subset of $\Omega$.
\end{lemma}

\proof Suppose the assertion is false. There exists an invariant set
$A$ whose intersection with many horizontal fibers is neither the full
fiber nor empty, mod $m_\om$.  That is, for some $\eps>0$, the
$\Pi$-\me\ of
\begin{equation}
  \label{erg30}
  B_\eps := \rset{\om \in \Omega} {m_\om (A \cap I_\om) \in [\eps,
  1-\eps]}
\end{equation}
is positive. By the Poincar\'e Recurrence Theorem and the Lebesgue
Density Theorem it is possible to pick a $(s,\om) \in A \cap (I \times
B_\eps)$ that is recurrent to $I \times B_\eps$ and such that
$(s,\om)$ is a density point of $A \cap I_\om$ within $I_\om$. We
claim that there exist a sufficiently large $n$ and a multi-index
$\mathbf{i} \in \mathcal{I}^n$ for which
\begin{equation}
  \label{erg35}
  m_\om( A \cap I_{\om,\mathbf{i}} ) > (1-\eps) \, m_\om(
  I_{\om,\mathbf{i}} )
\end{equation}
and 
\begin{equation}
  \label{erg36}
  \ma^n ( I_{\om,\mathbf{i}} ) = I_{\om_n} \subset I \times B_\eps. 
\end{equation}
In fact, among the infinitely many $n$ that verify $\ma^n(s,\om) \in I
\times B_\eps$, we can choose, by Lemma \ref{lem-hyp}, one for which
$I_{\om,\mathbf{i}_n (s,\om)}$ is so small that (\ref{erg35}) is
verified for $\mathbf{i} = \mathbf{i}_n (s,\om)$. The equality in
(\ref{erg36}) is true by the Markov property of $\phi(\cdot,
\om_{n-1}) \circ \cdots \circ \phi(\cdot, \om)$ (having used again
notation (\ref{sk-omk})).

Since the restriction of $\ma^n$ to $I_{\om,\mathbf{i}}$ is linear and
$A$ is invariant, we deduce from (\ref{erg36}) that $m_{\om_n} (A \cap
I_{\om_n}) > 1-\eps$, which contradicts (\ref{erg30}), because $\om_n
\in B_\eps$.  Therefore an invariant set, mod $\mu$, can only occur in
the form $I \times B$. That $B$ is measurable is a consequence of the
next lemma.  
\qed

\begin{lemma}
  For $i=1,2$, let $(\Sigma_i, \mathscr{A}_i, \nu_i)$ be two \pr\
  spaces, the second of which complete. If $B_1 \in \mathscr{A}_1$,
  $\nu_1(B_1) > 0$, and $B_1 \times B_2 \in \mathscr{A}_1 \otimes
  \mathscr{A}_2$, then $B_2 \in \mathscr{A}_2$.
\end{lemma}

\proof See \cite[Lemma A.1]{l1}.

\begin{theorem}
  \label{thm-erg}
  $(\ps, \mu, \ma)$ is \erg.
\end{theorem}

\proof Suppose the \sy\ is not \erg. By Lemma \ref{lem-erg}, we have
an invariant set $I \times B$, with $\Pi(B) \in (0,1)$. Set $B^c :=
\Omega \setminus B$. Without loss of generality, we assume that
the \tr y of \emph{no} point of $I \times B$ intersects $I \times B^c$
(otherwise it is easy to modify the following argument to deal with a
negiglible set of exceptions).

By the \erg ity of $(\Omega, \Pi, \{ \tau_z \})$ and the almost sure
irreducibility of the \ra\ \en, one can find an $\om \in B$ which is
transitive in the sense of (\ref{as-irr}), and a $y \in \Z^d$, such
that
\begin{equation}
  \label{erg40}
  \tau_y \, \om \in B^c.
\end{equation}
The transitivity of $\om$ means that there exist $n \in \Z^+$ and $s
\in I$ such that $X_n(s,\omega) = y$, whence, using (\ref{sk-omk}),
\begin{eqnarray}
  \ma^n(s,\om) &=& (s_n, \tau_{D(s_{n-1},\om_{n-1})} \circ \cdots 
  \circ \tau_{D(s,\om)} \, \om) \nonumber \\
  &=& (s_n, \tau_{X_n(s,\om)} \, \om) \nonumber \\
  &=& (s_n, \tau_y \, \om) \in I \times B^c,
\end{eqnarray}
the inclusion descending from by (\ref{erg40}). This contradicts the
initial assumption.  
\qed

We now prove a strong stochastic property for the specific
(vector-valued) observable $D$, namely the CLT for the family $\{ D
\circ \ma^k \}$. In the language of RWREs, this can be formulated as
follows.

\begin{theorem}
  \label{thm-clt}
  The stochastic process $\{ X_n \}$, defined in
  \emph{(\ref{def-xn-dk})}, satisfies the CLT with mean zero and
  covariance matrix $C := \{ c_{\alpha \beta} \}_{\alpha,\beta =
    1}^d$, where
  \begin{displaymath}
    c_{\alpha \beta} := \int_\ps D^{(\alpha)}(s,\om)
    D^{(\beta)}(s,\om) \: \mu (ds d\om),
  \end{displaymath}
  having denoted by $D^{(\alpha)}$ the $\alpha$-th component of the
  vector $D$.
\end{theorem}

\proof By elementary martingale theory \cite{hh} it suffices to prove
that $\{ X_n \}$ is a (multidimensional) martingale whose increments
$D_n = X_{n+1} - X_n$ have covariance matrix $C$ for all $n$ (if a
specific reference is needed, the first theorem of \cite{w} implies
the result).

From the considerations outlined in the beginning of Section
\ref{sec-erg-clt} it is not hard to see that $\mathscr{F}_n$, the
$\sigma$-algebra generated by $X_1, \ldots, X_n$ (equivalently, by
$D_0, \ldots, D_{n-1}$) corresponds to the partition $\mathscr{D}_n :=
\{ A_\mathbf{i} \}_{\mathbf{i} \in \mathcal{I}^n}$, with
\begin{equation}
  A_\mathbf{i} := \bigcup_{\om \in \Omega} I_{\om, \mathbf{i}} \, .
\end{equation}
Therefore, denoting by $\mathbb{E}_\mu$ the conditional expectation
w.r.t.\ $\mu$, we have
\begin{eqnarray}
  \mathbb{E}_\mu (D_n | \mathscr{F}_n) &=& \sum_{\mathbf{i} \in
  \mathcal{I}^n} \left[ \frac1 {\mu (A_\mathbf{i})} \int_{A_\mathbf{i}}
  D_n (s,\om) \, \mu( ds d\om ) \right] 1_{A_\mathbf{i}} = \nonumber \\
  \label{clt10}
  &=& \sum_{\mathbf{i} \in \mathcal{I}^n} \left[ \frac1 {\mu
  (A_\mathbf{i})} \int_\Omega \int_{I_\mathbf{i} (\om)} D_n (s,\om) \,
  ds \, \Pi(d\om) \right] 1_{A_\mathbf{i}} ,
\end{eqnarray}
where $1_{A_\mathbf{i}}$ is the indicator \fn\ of $A_\mathbf{i}$.
For $\mathbf{i}, \om$ fixed and $s$ ranging in $I_{\mathbf{i}}$ (we
are dropping the dependence on $\om$ from all the notation), the first
$n$ positions of the walk are determined, say by the values $x_k :=
X_k(s,\om)$ $(k=0, 1, \ldots, n)$. Thus, the inner integral in
(\ref{clt10}) becomes
\begin{eqnarray}
  \int_{I_\mathbf{i}} D_n (s,\om) \, ds &=& \sum_{j=1}^N d_j \, m \!
  \left( I_{(\mathbf{i}, j)} \right) \nonumber \\ &=& \sum_{j=1}^N d_j
  \, p_{0x_1} \cdots\, p_{x_{n-1}, x_n} \, p_{x_n, x_n+d_j} = 0,
\end{eqnarray}
by (\ref{lambda})-(\ref{zero-l-drift}). Hence the r.h.s.\ of
(\ref{clt10}) is identically zero, proving the martingale property. As
concerns the covariances of $D_n = D \circ \ma^n$, their constancy in
$n$ follows from the invariance of $\mu$.  
\qed

\section{Cocycles and recurrence}
\label{sec-rec}

The upcoming definitions and results apply to a general \dsy.

\begin{definition}
  \label{def-co}
  Let $(\Sigma, \nu, F)$ be a \pr-preserving dynamical \sy, and $f$ a
  measurable \fn\ $\Sigma \longrightarrow \Z^d$. The family of \fn s
  $\{ S_n \}_{n\in \N}$, defined by $S_0(\xi) \equiv 0$ and, for $n \ge
  1$,
  \begin{displaymath}
    S_n(\xi) := \sum_{k=0}^{n-1} (f \circ F^k) (\xi)
  \end{displaymath}
  is called a \emph{commutative, $d$-dimensional, discrete cocycle}
  or, more precisely, the \emph{cocycle of $f$}.
\end{definition}

\begin{definition}
  \label{def-co-rec}
  The discrete cocycle $\{ S_n \}$ is called \emph{recurrent} if, for
  $\nu$-almost all $\xi \in \Sigma$, there exists a subsequence $\{
  n_j = n_j(\xi) \}$ such that
  \begin{displaymath}
    \forall j \in \N, \qquad S_{n_j}(\xi) = 0.
  \end{displaymath}
\end{definition}

A remarkable sufficient condition for cocycle recurrence was given by
Schmidt \cite{s}:

\begin{theorem}
  \label{thm-co-rec}
  Assume that $(\Sigma, \nu, F)$ is \erg\ and denote by $R_n$ the
  distribution of $S_n / n^{1/d}$, relative to $\nu$.  If there exists
  a positive-density sequence $\{ n_k \}_{k\in\N}$ and a constant
  $K>0$ such that
  \begin{displaymath}
    R_{n_k}( \mathcal{B}(0,\rho) ) \ge K \rho^d 
  \end{displaymath}
  for all sufficiently small balls $\mathcal{B}(0,\rho) \subset \R^d$
  (of center 0 and radius $\rho$), then the cocycle $\{ S_n \}$ is
  recurrent.
\end{theorem}

\begin{remark}
  Schmidt proved the above only in the case where $F$ is an
  automorphism (i.e., it is invertible mod $\nu$) \cite{s}. It is easy,
  however, to extend his proof to the full generality claimed by
  Theorem \ref{thm-co-rec}. See Section \ref{subs-ne} of the Appendix.
\end{remark}

In dimension 1 and 2, if $\{ S_n \}$ satisfies the CLT with zero mean
(even with a degenerate covariance matrix), it satisfies the main
hypothesis of Theorem \ref{thm-co-rec}. Therefore, coming back to our
\sy, since $\{ X_n \}$ is the cocycle of $D$ via (\ref{def-xn-dk}),
and in view of Theorems \ref{thm-erg} and \ref{thm-clt}, we obtain the
following

\begin{theorem}
  \label{thm-rec}
  If $d \le 2$, the RWRE described in Section \ref{sec-intro} is
  \emph{almost surely recurrent}. This means that, for $\Pi$-a.e.\
  $\om$, the \rw\ $\{ X_n \}$, subject to the law $P_\om$, verifies
  $X_{n_j} = 0$, with \pr\ 1, for a subsequence $\{ n_j \}$.
\end{theorem}

\section{Conclusions and generalizations}
\label{sec-gen}

We have presented a fairly natural way---at least to a \hyp\
dynamicist---to represent a RWRE as a \pr-preserving \dsy. This is
achieved by implementing the local dynamics of the particle in terms
of one-dimensional perfect Markov maps and then considering a sort of
``union'' of all of them. (Notice, however, that the \sy\ we have
introduced in Section \ref{sec-povp} is not the only one that is
called `the point of view of the particle' in the field.)

In this representation, elementary considerations of \erg\ theory can
produce interesting, new and old, results in a nearly effortless way:
e.g., in the martingale case, the CLT with no ellipticity assumption
and, using a powerful theorem by Schmidt, recurrence in dimension $d
\le 2$.

In fact, the results contained in this note can be generalized in a
number of ways, from the more to the less evident. For example:

\begin{enumerate}
\item There is no need for $p_{xy}$ to be zero when $y-x \not\in
  \Lambda$, as long as it decays so fast that the \fn\ $D$ of
  (\ref{def-d}) has square-integrable modulus. The condition $\sum_y
  |y|^2 \, p_{0y} \le K$, for a.a.\ $\om$, suffices.

\item If we know that the \ra\ \en\ is \erg\ for the action of a
  subgroup $\Gamma \subset \Z^d$ (i.e., $(\Omega, \Pi, \{ \tau_z \}_{z
  \in \Gamma})$ is \erg, which is a stronger hypothesis than the one
  we made in Section \ref{sec-intro}), we can relax the condition of
  almost sure irreducibility of $p$ and ask that (\ref{as-irr}) be
  verified only for $y \in \Gamma$. (In fact, the proof of Theorem
  \ref{thm-erg} works equally well if $\{ \tau_z \}$ is restricted to
  $z \in \Gamma$.) In particular, if $(\Omega, \Pi)$ is an i.i.d.\
  \en, it suffices to require the existence of a (one-dimensional)
  subgroup $\Gamma$ such that (\ref{as-irr}) holds $\forall y \in
  \Gamma$.

\item More importantly, the bistochasticity condition (\ref{bi-en})
  can be done away with. 

\item Finally, one can prove not only the annealed CLT, but the
  \emph{quenched} invariant principle, i.e., for a \emph{fixed}
  typical \en\ $\om$, the convergence of the rescaled \tr y
  \begin{equation}
    t \mapsto \frac1{\sqrt{n}} \sum_{k=1}^{[nt]} X_k
  \end {equation}
  to a $d$-dimensional Brownian motion, for $t \in [0,1]$. The
  annealed invariant principle follows.
\end{enumerate}

These advances will be presented in \cite{l3}.
\skippar

The techniques exposed in this paper can also be used separately. For
instance, one can apply Theorem \ref{thm-co-rec} to the \rw s of
\cite{ko}. These are more general bistochastic RWREs than the ones
discussed here, in that they do not have zero local drift as in
(\ref{zero-l-drift}), but zero \emph{mean} drift:
\begin{equation}
  \label{zero-drift}
  \int_\Omega \sum_y y\, p_{0y}(\om) \: \Pi(d\om) = 0.
\end{equation}
Subject to a uniform ellipticity condition, namely that, $\Pi$-almost
surely,
\begin{equation}
  \label{ko-ell}
  p_{0y} \ge \eps, \qquad \forall y \in \Lambda;
\end{equation}
and to an additional, rather cumbersome but essential, hypothesis they
call \emph{condition (H)}, Komorowski and Olla \cite[Thm.~2.2]{ko}
prove the annealed zero-mean CLT for $\{ X_n \}$, i.e., they prove
Theorem \ref{thm-clt} for the \dsy\ $(\ps, \mu, \ma)$ adapted to their
case. Since the \erg ity of that \sy\ is guaranteed by (\ref{ko-ell})
(as in Theorem \ref{thm-erg}), Theorem \ref{thm-rec}, that is the a.s.\
recurrence in dimension 1 and 2, holds true in this case as well.

\appendix

\section{Appendix}

\subsection{An example of a bistochastic environment with zero local drift}
\label{subs-ex}

As it may not be immediately intuitive how to construct a \ra\ \en\
satisfying the assumptions of Section \ref{sec-intro}, we give here an
explicit example, emphasizing that the $(\Omega, \Pi)$ produced below
is not as general as we are able to deal with in the present paper
(see also Section \ref{sec-gen}).

We start by fixing a finite $\Lambda_0 \subset \Z^d$, of cardinality
$N_0$.  Calling $\mathbb{B} = \mathbb{B} (\Lambda_0)$ the set of all
bistochastic matrices indexed by the elements of $\Lambda_0$, it is
known that $\mathbb{B}$ can be identified with a convex polytope of
$\R^{(N_0-1)^2}$ (each matrix has $N_0^2$ entries and $2N_0 - 1$
independent conditions on them). The Birkhoff--Von Neumann Theorem
states that any element of $\mathbb{B}$ can be expressed as a convex
combination of the $N_0!$ permutation matrices of $\mathbb{B}$, which
are thus the extremal points of the polytope (the permutation matrices
are those that are obtained by permuting the columns of the identity
matrix) \cite{b}.

Now endow $\mathbb{B}$ with any absolutely continuous \pr\ $\pi_0$
(w.r.t.\ the Lebesgue \me\ in $\R^{(N_0-1)^2}$) and define $(\Omega,
\Pi) := (\mathbb{B}, \pi_0)^{\Z^d}$, in the sense of the tensor
product of \me\ spaces. The generic element of $\Omega$ is denoted
$\om = \{ \om^{(\zeta)} \}_{\zeta \in \Z^d}$, where $\om^{(\zeta)} =
\{ \om^{(\zeta)}_{xy} \}_{x, y \in \Lambda_0} \in \mathbb{B}$. To keep
the notation simple, pretend that $\om^{(\zeta)}_{xy}$ exists for all
$x,y \in \Z^d$, equaling 0 where not otherwise defined.

Setting
\begin{equation}
  \label{re10}
  o_{xy} = o_{xy}(\om) := \frac1{N_0} \sum_{\zeta \in \Z^d} 
  \om^{(\zeta)}_{x+\zeta, y+\zeta}
\end{equation}
gives rise to a bistochastic \en\ $o = o(\om) := \{ o_{xy} \}_{xy}$,
as it can easily be verified. If we further set 
\begin{equation}
  \label{re20}
  p_{xy} = p_{xy}(\om) := \frac{o_{xy} + o_{x,-y}} 2
\end{equation}
and $p = p(\om) := \{ p_{xy} \}_{xy}$, we obtain a bistochastic \en\
which is also \emph{balanced} in the sense of Lawler \cite{la}, namely
$p_{xy} = p_{x,-y}$. This is a special case of the zero-local-drift
condition (\ref{zero-l-drift}). Also, (\ref{lambda}) holds true for
$\Lambda := \Lambda_0 \cup -\Lambda_0$.

For $z \in \Z^d$, we define $\tau_z: \Omega \longrightarrow \Omega$
via
\begin{equation}
  \label{re30}
  (\tau_z \om)^{(\zeta)} := \om^{(\zeta-z)}.
\end{equation}
Since the $\om^{(\zeta)}$ are i.i.d.\ \ra\ matrices, obviously $\{
\tau_z \}$ leaves $\Pi$ invariant and is \erg. Equality
(\ref{zd-action}) is easily checked for $o$ and thus for $p$.

Finally, since $\pi_0$ has a density on $\mathbb{B} \subset
\R^{(N_0-1)^2}$, the \pr\ that $\om^{(\zeta)}$ has a zero entry is
null, which proves the (nonuniform) ellipticity of $o$ and $p$,
implying (\ref{as-irr}).

\subsection{Partial proof of Schmidt's Theorem}
\label{subs-ne}

In this section we show that Theorem \ref{thm-co-rec}, which was
proved by Schmidt only for $F$ invertible \cite{s}, easily extends to
the case of a general endomorphism $F$ (note that Conze has
independently given a weaker result than Theorem \ref{thm-co-rec}
which does not require invertibility \cite{c}).

In fact, if $(\Sigma, \mathscr{A}, \nu, F)$ is a \pr-preserving,
noninvertible, \dsy\ (here $\mathscr{A}$ is the $\sigma$-algebra
defined on $\Sigma$), one can consider its \emph{natural extension} in
the sense or Rohlin \cite{r}. Shying away from the details of its
construction (which can be found, e.g., in \cite{r, cfs}), we simply
recall that the natural extension of $(\Sigma, \mathscr{A}, \nu, F)$
is a \pr-preserving invertible \dsy\ $(\bar{\Sigma},
\bar{\mathscr{A}}, \bar{\nu}, \bar{F})$ for which there exists a
measurable projection $\pi : \bar{\Sigma} \longrightarrow \Sigma$ with
the following properties: First, the commutation condition
\begin{equation}
  \label{ne10}
  \pi \circ \bar{F} = F \circ \pi,
\end{equation}
which explains in what sense $\bar{F}$ extends $F$. Second, setting
\begin{equation}
  \label{ne20}
  \bar{\mathscr{A}}_0 := \pi^* \mathscr{A} := \lset{\pi^{-1} A} 
  {A \in \mathscr{A}},
\end{equation}
one has that $(\bar{\Sigma}, \bar{\mathscr{A}}_0, \bar{\nu})$ and
$(\Sigma, \mathscr{A}, \nu)$ are isomorphic \me\ spaces. Furthermore
$\bar{\mathscr{A}}_0$ is $\bar{F}$-invariant (this means, using again
the terminology of (\ref{ne20}), that $\bar{F}^* \bar{\mathscr{A}}_0
\subset \bar{\mathscr{A}}_0$) but in general not
$\bar{F}^{-1}$-invariant.

Hence any $\mathscr{A}$-measurable \fn\ $f: \Sigma \longrightarrow \R$
is isomorphically associated to the $\bar{\mathscr{A}}_0$-measurable
\fn\ $\bar{f} := f \circ \pi: \bar{\Sigma} \longrightarrow \R$.
Moreover, $\forall k \in \N$, $f \circ F^k$ and $\bar{f} \circ
\bar{F}^k$ have the same distribution (because, by (\ref{ne10}),
$\overline{f \circ F^k} = \bar{f} \circ \bar{F}^k$). Finally, it is
known that $(\bar{\Sigma}, \bar{\mathscr{A}}, \bar{\nu}, \bar{F})$ is
\erg\ \iff $(\Sigma, \mathscr{A}, \nu, F)$ is \erg\ \cite{cfs}.

Coming back to Theorem \ref{thm-co-rec} for a noninvertible $F$, one
can apply Schmidt's proof to the natural extension, which is \erg, and
its cocycle $\{ \bar{S}_n \}$ (obvious definition). The latter
verifies the main hypothesis of the theorem because $\{ S_n \}$ does,
as explained above. The fact that $\{ S_n \}$ is recurrent \iff $\{
\bar{S}_n \}$ is concludes the proof.
\qed

\end{document}